\newif
\ifconfver

\confvertrue        

\ifconfver
    \documentclass[10pt,journal,final,twoside]{IEEEtran}
\else
    \documentclass[11pt,draftcls,onecolumn]{IEEEtran}
\fi

\usepackage[utf8]{inputenc}
\usepackage{xspace,empheq,fancybox,amsmath,amssymb,amsfonts,epstopdf,epsfig,syntonly,times,amsthm,euscript} 
\usepackage{psfrag,color,bm,array,tabularx, subfigure}
\usepackage{cite,url,footnote,xspace,syntonly,algorithmic}
\usepackage{verbatim,multirow}

\usepackage[linesnumbered,ruled,vlined]{algorithm2e}
\usepackage[T1]{fontenc}

\usepackage{textcomp}
\usepackage{xcolor}
\usepackage{mathtools}
\usepackage{makecell}
\usepackage{bbm, booktabs, makecell, subcaption, graphicx, titlesec}
\usepackage{microtype} 
\newcolumntype{M}[1]{>{\centering\arraybackslash}m{#1}}
\newcolumntype
{N}{@{}m{0pt}@{}}
\DeclareMathOperator*{\argmax}{\arg\max}
\def\BibTeX{{\rm B\kern-.05em{\sc i\kern-.025em b}\kern-.08em
    T\kern-.1667em\lower.7ex\hbox{E}\kern-.125emX}}
\newcommand{\KB}{\color{black}{}}    
\newcommand{\KBB}{\color{black}{}}  
\newcommand{\KBA}{\color{black}{}}
\newcommand{\FINAL}{\color{black}{}}

\input{mysymbol.sty}
\allowdisplaybreaks

\IEEEoverridecommandlockouts
\title{\LARGE \bf Risk-Constrained Reinforcement Learning \\ for Inverter-Dominated Power System Controls}
\author{\IEEEauthorblockN{Kyung-bin Kwon, Sayak Mukherjee, Thanh Long Vu, and Hao Zhu}

\thanks{K. Kwon and H. Zhu are with the Chandra Department of Electrical and Computer
Engineering, The University of Texas at Austin, 2501 Speedway, Austin,
TX, 78712, USA; Emails: \{kwon8908kr, haozhu\}@utexas.edu.}
\thanks{S. Mukherjee and T.~L. Vu are with the Optimization and Control Group, Pacific Northwest National Laboratory, Richland, WA 99352, USA, Emails: \{sayak.mukherjee, thanhlong.vu\}@pnnl.gov.}%
\thanks{K. Kwon was a PhD intern with the Optimization and Control Group, Pacific Northwest National Laboratory.}%
\thanks{\protect\rule{0pt}{3mm} K. Kwon and H. Zhu have been supported by NSF Grants 2130706, 2150571, and by ARO Grant W911NF2310266. The research of S. Mukherjee and T.~L. Vu  is part of the E-COMP (Energy System Co-Design with Multiple Objectives and Power Electronics) Initiative conducted under the Laboratory Directed Research and Development program at Pacific Northwest National Laboratory, a multiprogram national laboratory operated by Battelle for the U.S. Department of Energy. }
}

\begin{document}
\maketitle

%

\begin{abstract}
This paper develops a risk-aware controller for grid-forming inverters (GFMs) to minimize large frequency oscillations in GFM inverter-dominated power systems. To tackle the high variability from loads/renewables, we incorporate a mean-variance risk constraint into the classical linear quadratic regulator (LQR) formulation for this problem. The risk constraint aims to bound the time-averaged cost of state variability and thus can improve the worst-case performance for large disturbances. The resulting risk-constrained LQR problem is solved through the dual reformulation to a minimax problem, by using a reinforcement learning (RL) method termed as stochastic gradient-descent with max-oracle (SGDmax). In particular, the zero-order policy gradient (ZOPG) approach is used to simplify the gradient estimation using simulated system trajectories. Numerical tests conducted on the IEEE 68-bus system have validated the convergence of our proposed SGDmax for GFM model and corroborate the effectiveness of the risk constraint in improving the worst-case performance while reducing the variability of the overall control cost.
\end{abstract}

\begin{IEEEkeywords}
Frequency control, grid-forming inverter (GFM), inter-area oscillations, mean-variance risk constraint, reinforcement learning (RL)
\end{IEEEkeywords}

%

\section{Introduction} \label{sec:IN}

Grid-forming inverters (GFMs) {\KBA are} increasingly important for establishing grid voltage and frequency in next-generation power systems with high penetration of low-carbon energy resources \cite{anttila}.  
Photovoltaics, wind generators, and energy storage devices lack in conventional primary and secondary controls as synchronous generators (SGs), and thus their integration greatly challenges grid stability. Advanced GFM technology can address this issue as they operate as independent voltage sources to support grid stability by controlling the voltage and frequency at the interfaces of new resources \cite{long}. 


{\KB The existing GFM control strategies consist of three main categories: droop control \cite{rath, wei, peyghami}, virtual synchronous generators \cite{ebrahimi, serban}, and virtual oscillator control \cite{grob}.} Especially, droop control is a well-established method to mitigate voltage and frequency fluctuations by following $P$-$\omega$ and $Q$-$V$ droop curves. By observing active and reactive powers from the network as inputs, it can vary the terminal frequency and voltage  depending on the internal voltage/power set-points \cite{du}. Thus, changing these set-points can affect the overall grid dynamics to quickly attain the steady-state operations after huge external perturbations due to, {\KBA e.g.}, sudden changes of the load/renewable. As a local control design, it is known that the overall  performance of multiple droop-controllers could  degrade for reducing inter-area oscillations in large-scale interconnection \cite{rath}. A decentralized control design among all GFMs can address this issue, with  state information exchange among GFMs  as enabled by communication network \cite{behrooz}. The number of communication links are typically limited, and thus a structured feedback design per the information-exchange graph among GFMs will be adopted later on.

Recent advances in data-driven methods, including both model-based and model-free ones, have provided significant advantages for solving optimal control problems. To design the decentralized GFM controller, there have been several data-driven techniques such as reinforcement learning (RL) \cite{li} and adaptive dynamic programming \cite{mukherjee2021scalable}. While these model-free approaches do not require to know the system's mathematical model,  they are known to suffer from the sample complexity issue which needs extensive data samples and large training time\cite{lintao}. {\KBA Thus, we choose a model-based approach by simulating the underlying system dynamics, which will greatly accelerate the policy search in practice \cite{mbrl1, mbrl2}.}

Moreover, a linear quadratic regulator (LQR) objective \cite{eskandari} typically used for the GFM control problem could unfortunately face severe performance degradation particularly in the worst-case oscillations. This is because it solely focuses on reducing the expected cost over time, failing to account for the highest oscillation scenarios. The latter issue is particularly pronounced in large-scale interconnections, where sudden load/renewable changes can trigger the inter-area oscillations that are known to experience high frequency fluctuations and poor damping levels. {\FINAL While some existing work \cite{eskandari} has attempted to use model-based RL for GFM controls, it is still an open question to account for the risk associated with extreme scenarios.}

To this end, we develop a risk-aware RL approach for the GFM control design problem while aiming to address large perturbations in GFM inverter-dominated systems. We formulate it as a constrained LQR problem with the so-termed mean-variance risk constraint. The latter is imposed on the overall deviations of state cost from its expectation, which can reduce the level of high system variability as a result of significant disturbances to enhance the worst-case performance. {\KBA To solve this problem, we implement an RL-based algorithm termed as stochastic gradient-descent with max-oracle (SGDmax), which utilizes zero-order policy gradient (ZOPG) as estimated gradients for reduced computational complexity \cite{kwon}. Expanding upon \cite{kwon} that utilize SGDmax in load frequency control (LFC) with simple dynamics in small microgrids system, we develop the system dynamics model by integrating detailed SG and GFM dynamics. In addition, we demonstrate the effectiveness of the method based on the widely recognized IEEE 68-bus system.}

{\KBA Our main contributions are three-fold. First, we formulate the GFM problem as an LQR problem with linearized dynamics integrating both SGs and GFMs through network coupling, {\KBB which is essential for incorporating the risk constraint}. Second, we consider high load perturbation into system dynamics that leads to increased system variability. To address this challenge, we incorporate the mean-variance risk constraint into LQR formulation and solve it using the RL-based SGDmax algorithm. 
Third, we demonstrate the effectiveness of the proposed method through numerical tests in the presence of high external disturbances. Specifically, we train an RL policy by generating the system trajectory offline and test it in the actual nonlinear system to demonstrate the improved efficacy of our proposed method.}

\section{System Modeling} \label{sec:SM}

We consider a power system consisting of $N_a$ areas with a total of $N_g$ synchronous generators (SGs) and $N_f$ grid-forming inverters (GFMs).
We denote SGs and GFMs as indexed by the sets $\mathcal{G} = \{1,2,\ldots, N_g\}$ and $\mathcal{F} = \{N_g+1, N_g+2, \ldots, N_g+N_f\}$, respectively. Without loss of generality (Wlog), every load bus is assumed to be connected to one GFM, as other load buses can be reduced. 


To model the overall system, we first present the dynamics of each SG and GFM here. The dynamics of SG $i \in \mathcal{G}$ is represented by the following swing equations \cite{kundur}:
\begin{subequations}
\begin{align}
&\dot{\delta}_i = \omega_i - \omega_0,\\
&\dot{\omega}_i = \frac{1}{M_i}\big[D_i(\omega_0-\omega_i) + P_i - P^n_i\big]
\end{align} \label{eq:sg}
\end{subequations}
\noindent where parameters $\{M_i, D_i\}$ are the inertia and damping coefficients, respectively. Let the vector $\mathbf{x}_i = [\delta_i, \omega_i]$ collect the  internal state of rotor angle and speed for SG $i$. {\FINAL The active power delivered to the network $P^n_i$ is compared with its generated one $P_i$ to  determine the rate of frequency deviation.}
Here, we do not incorporate the frequency regulation of the automatic generation control (AGC) due to its slow time-scale. However, the AGC extension is possible to include as part of the network power flow as considered later on. 

\begin{figure}[t]
	\centering
	\includegraphics[width=\linewidth]{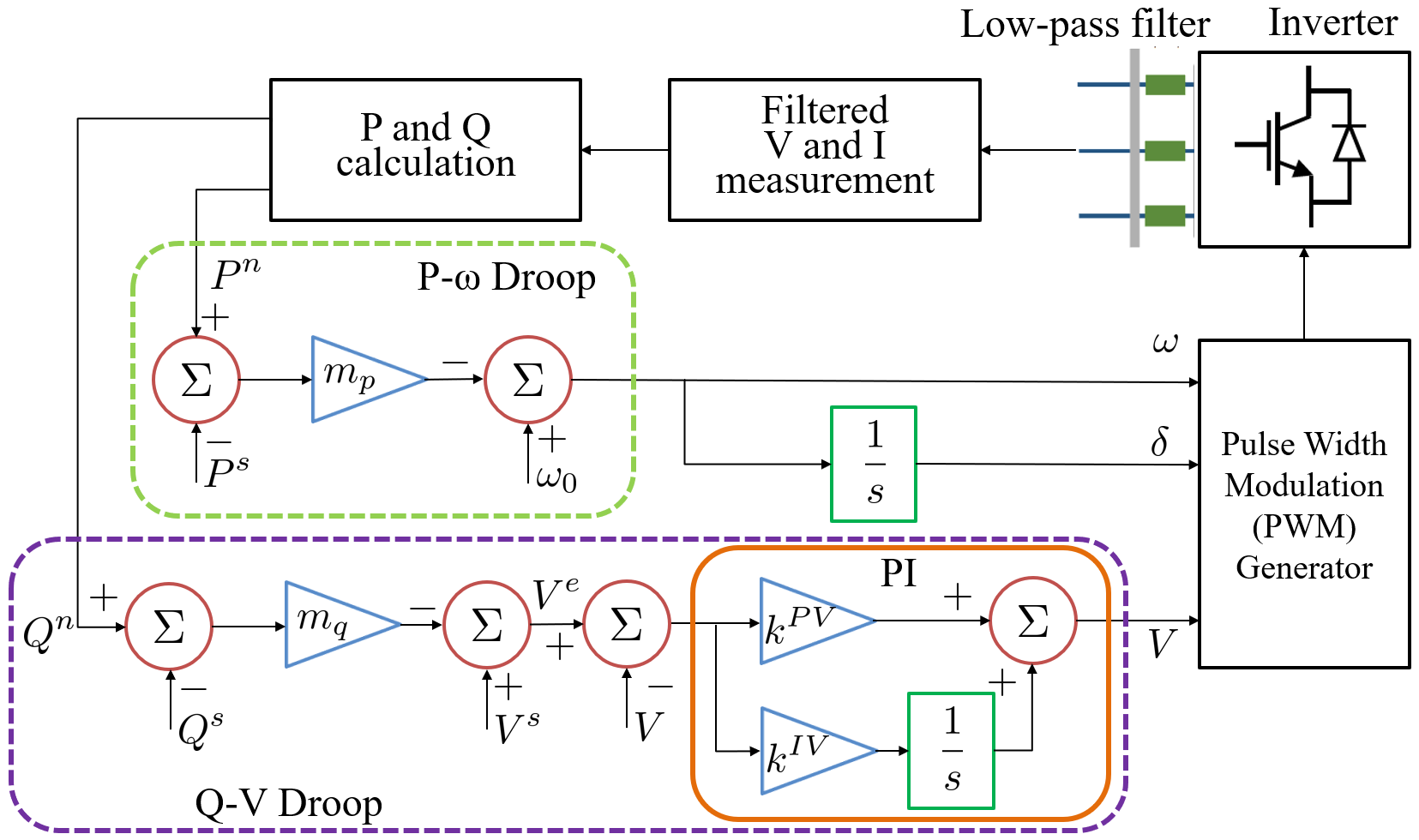}
	\caption{GFM dynamics based on droop controls.}
	\label{fig:gfm}
\end{figure}

Each of the $N_f$ GFMs acts as a controllable voltage source at the DER-connected load bus \cite{du}. As illustrated in Fig.~\ref{fig:gfm}, the internal dynamics  utilizes the  $P$-$\omega$ and $Q$-$V$ droop control curves \cite{behrooz}.
For each GFM $j\in\mathcal{F}$, the active and reactive powers delivered by the network, namely $P_j^n$ and $Q_j^n$, are calculated using the terminal voltage/current measurements that go through a low-pass filter. The droop controller then uses the difference from the corresponding power set-point to determine the actuation signal. For example, the difference between $P_j^n$ and the active power set-point $P^s_j$ is multiplied by the droop gain $m^p_j$ to determine the signal $\omega_j$; and similarly for the voltage error signal $V^e_j$. 
Note that the $Q$-$V$ droop is also followed by a proportional-integral (PI) controller to further regulate the deviations of the voltage error $V^e_j$, with $k^{pv}_j$ and $k^{iv}_j$ as the proportional and integral gains, respectively. Last, $\omega, \delta$ and $V$ are sent to the pulse width modulation (PWM) generator and the frequency and voltage of the node are set accordingly by the inverter.
Thus, the dynamic model of GFM $j \in \mathcal{F}$ can be expressed as follows:
\begin{subequations}
\begin{align}
&\dot{\delta}_j = \omega_j - \omega_0,\\
&\dot{\omega}_j = \frac{1}{\tau_j}\big[\omega_0 - \omega_j +m^p_j (P^s_j - P^n_j)\big],\\
&\dot{V}^e_j = \frac{1}{\tau_j}\big[V^s_j - V^e_j - V_j + m^q_j (Q^s_j - Q^n_j)\big],\\
&\dot{V}_j = k^{pv}_j \dot{V}^e_j + k^{iv}_j V^e_j
\end{align} \label{eq:gfm}
\end{subequations}
\noindent where $\tau_j$ is a pre-determined droop time constant. The state vector per GFM $j$ becomes $\mathbf{x}_j := [\delta_j, \omega_j, V^e_j, V_j]$. Thus, the GFM works by controlling its $V_j$ and $\omega_j$ via adjusting the voltage and power set-points, which are included by the vector $\mathbf{u}_j := [V^s_j, P^s_j. Q^s_j]$.

\setlength{\textfloatsep}{0.35cm}

Based on \eqref{eq:sg} and \eqref{eq:gfm}, the dynamics of SGs and GFMs are coupled by the network power flow (PF) to determine $\{P_\ell^n\}$ and $\{Q_\ell^n\}$ for $\ell\in \mathcal{G} \cup \mathcal{F}$. Thus, we can establish the overall system dynamics through the steady-state PF analysis. By employing the Kron reduction \cite{kron}, we first eliminate all other buses except for the SG and GFM ones, to consider the PF among $\{V_\ell, \delta_\ell\}_{\ell\in\mathcal{G}\cup\mathcal{F}}$ as given by:
\begin{align}
P^n_\ell &= \sum_{k=1}^{N_g+N_f} V_\ell V_k(G_{\ell k}\cos(\delta_\ell - \delta_k)+B_{\ell k}\sin(\delta_\ell - \delta_k))\nonumber\\
Q^n_\ell &= \sum_{k=1}^{N_g+N_f} V_\ell V_k(G_{\ell k}\sin(\delta_\ell - \delta_k)-B_{\ell k}\cos(\delta_\ell - \delta_k))\nonumber
\end{align}

{\KBB To formulate the linearized dynamics required for the risk constraint reformulation, we linearize the PF equations around the steady-state operating point, resulting in}
\begin{align}
    \begin{bmatrix} \Delta \mathbf{P}^g \\ \Delta \mathbf{P}^f \\ \Delta \mathbf{Q}^f \end{bmatrix} &= \begin{bmatrix} 
     \frac{\partial \mathbf{P}^g}{\partial \bm{\delta}^g} & \frac{\partial \mathbf{P}^g}{\partial \bm{\delta}^f} & \frac{\partial \mathbf{P}^g}{\partial  \mathbf{V}^f} 
     \\ \frac{\partial \mathbf{P}^f}{\partial \bm{\delta}^g} & \frac{\partial \mathbf{P}^f}{\partial \bm{\delta}^f} & \frac{\partial \mathbf{P}^f}{\partial \mathbf{V}^f}
     \\ \frac{\partial \mathbf{Q}^f}{\partial \bm{\delta}^g} & \frac{\partial \mathbf{Q}^f}{\partial \bm{\delta}^f} & \frac{\partial \mathbf{Q}^f}{\partial \mathbf{V}^f} \end{bmatrix} 
     \begin{bmatrix} \Delta\bm{\delta}^g\\ \Delta\bm{\delta}^f\\ \Delta \mathbf{V}^f \end{bmatrix} \label{eq:lin}
\end{align}
\noindent where bold notation indicates the vector form by concatenating all corresponding scalar variables, with $g$ and $f$ indicating the SG and GFM components, respectively. Note that the reactive power component from SGs is not considered as we do not model the SGs' exciter control. 

By substituting \eqref{eq:lin} into \eqref{eq:sg} and \eqref{eq:gfm}, we can formulate the overall dynamics in continuous-time as
\begin{align}
    \dot{\mathbf{x}} = \mathbf{A}_c \mathbf{x} + \mathbf{B}_c \mathbf{u} + \bm{\xi} \label{eq:dynamics1}
\end{align}
\noindent where  $\mathbf{x}\!:=\![\Delta {\bm{\delta}}_g,\!\Delta {\bm{\omega}}_g,\!\Delta {\bm{\delta}}_f,\!\Delta {\bm{\omega}}_f,\!\Delta {\mathbf{V}}^e_f,\!\Delta {\mathbf{V}}_f]^\intercal\!\in\!\mathbb{R}^{2N_g + 4N_f}$, and $\mathbf{u} := [\Delta \mathbf{V}^s_f, \Delta \mathbf{P}^s_f, \Delta \mathbf{Q}^s_f]^\intercal \in \mathbb{R}^{3N_f}$.
Due to the linearization, all the variables in $\mathbf{x}$ and $\mathbf{u}$ now represent the deviations from the corresponding steady-state values.
Note that we add $\bm{\xi}$ which denotes random perturbations to system states, such as external disturbance or imperfect system modeling.
By considering the GFM control time $\Delta t$, we represent the discrete-time dynamics based on $\eqref{eq:dynamics1}$ as
\begin{align}
\mathbf{x}_{t+1} = \mathbf{A} \mathbf{x}_t + \mathbf{B} \mathbf{u}_t + \bm{\xi}_t,  t = 0, 1, \ldots \label{eq:dynamics}
\end{align}
\noindent where $\mathbf{A} = \mathbf{I} + \Delta t \cdot \mathbf{A}_c$ and   $\mathbf{B} = \Delta t  \cdot \mathbf{B}_c$
with $\mathbf{I}$ and $\Delta t$ denoting the identity matrix and a small enough time step, respectively. {\KBA Note that we assume \eqref{eq:dynamics} to be stable with sufficient communication links in the network, which is reasonable since the system is linearized around the steady-state operating point.} 

%

\section{Risk-constrained GFM Problem} \label{sec:RC}

Under the system dynamics in \eqref{eq:dynamics}, we can formulate the GFM control problem as an optimal control one with the linear quadratic regulator (LQR) objective, by minimizing the total cost of state and control, as
\begin{align}
 &\min_{\mathbf{K} \in \mathcal{K}} \; R_0(\mathbf{K}) \! = \!\lim_{T \to \infty}\frac{1}{T}\mathbb{E} \sum_{t=0}^{T-1} \big[\mathbf{x}^\intercal_t \mathbf{Q} \mathbf{x}_t + \mathbf{u}^\intercal_t \mathbf{R} \mathbf{u}_t\big] \label{eq:obj}
 \end{align}
 \noindent where matrices $\{\mathbf{Q},\mathbf{R}\}$ are positive (semi-)definite matrices used to weight the state and control variables into a single cost. Our goal is to find the best structured controller gain matrix $\mathbf{K} \in \mathbb{R}^{3N_f \times (2N_g+4N_f)}$ that linearly maps from $\mathbf{x}_t$ to $\mathbf{u}_t$, namely $\mathbf{u}_t = - \mathbf{K} \mathbf{x}_t$. Here, $\mathcal{K}$ indicates the structured feedback set defined by the information-exchange graph, {\KBA leading to the design of a distributed control system based on information transmitted through communication links.} Specifically, for any GFM or SG node $\ell \in \mathcal{G} \cup \mathcal{F}$ and GFM node $j\in\mathcal{F}$, the structured set $\mathcal{K}$ is defined as
\begin{align}
    \mathcal{K} = \{\mathbf{K}: \mathbf{K}_{j,\ell} = 0 \;\text{if and only if}\; j \nleftrightarrow \ell) \}, \nonumber
\end{align}
\noindent with $j \nleftrightarrow \ell$ implying that nodes $\ell$ and $j$ are not connected through a communication link. Note that the size of $\mathbf{K}_{j,\ell}$ is different according to $\ell$, i.e. $\mathbf{K}_{j,\ell} \in \mathbb{R}^{3 \times 2}$ if $\ell \in \mathcal{G}$ and $\mathbf{K}_{j,\ell} \in \mathbb{R}^{3 \times 4}$ if $\ell \in \mathcal{F}$.
While the sparsity of $\mathbf{K}$ presents challenges in the analysis of the feasible region \cite{feng}, it will not affect the implementation of our proposed algorithm.

Although the LQR objective in \eqref{eq:obj} effectively reduces oscillations on average, focusing solely on the average trajectory performance cannot account for the substantial system variability.
Specifically, high load perturbations introduce additional fluctuations into the dynamics described in \eqref{eq:dynamics}, which we can redefine the perturbation term as $\xi'_t$ that includes both the initial noise $\xi_t$ and the fluctuations from significant load changes. 
This change results in an increase in the variability of the system trajectory, posing a significant challenge to the LQR-based design. Notably, this issue becomes critical in interconnected grids where inter-area oscillations can arise from these substantial disturbances, thereby reducing the worst-case damping performance of the controller.

To tackle this issue, we put forth a risk-constrained LQR formulation by limiting the so-termed mean-variance risk measure, as
\begin{align}
	&\min_{\mathbf{K} \in \mathcal{K}} \; R_0(\mathbf{K}) \! = \!\lim_{T \to \infty}\frac{1}{T}\mathbb{E} \sum_{t=0}^{T-1} \big[\mathbf{x}^\intercal_t \mathbf{Q} \mathbf{x}_t + \mathbf{u}^\intercal_t \mathbf{R} \mathbf{u}_t\big] \label{eq:opt}\\
	&\textrm{s.t.}~R_c(\mathbf{K}) \! =\!\! \lim_{T \to \infty}\!\frac{1}{T} \mathbb{E} \sum_{t=0}^{T-1} \big(\mathbf{x}_t^\intercal \mathbf{Q} \mathbf{x}_t - \mathbb{E}\big[\mathbf{x}_t^\intercal \mathbf{Q} \mathbf{x}_t \vert \mathcal{H}_t\big]\big)^2 \leq c. \nonumber
\end{align}
We denote $\mathcal{H}_t = [\mathbf{x}_0, \mathbf{u}_0, \ldots, \mathbf{x}_{t-1}, \mathbf{u}_{t-1}]$ as the system state and control trajectory up to time $t$, and $c$ as a risk tolerance parameter. Note that the risk constraint bounds the deviations of the realized state cost $\mathbf{x}_t^\intercal \mathbf{Q} \mathbf{x}_t$ from its expected value. This constraint enables us to mitigate the worst-case scenarios of very high system variability as caused by external load disturbances and imperfect modeling. {\KBB Here, we assume the feasibility of \eqref{eq:opt} by assuming the system stability and sufficiently large threshold $c$ with adequate connectivity.}


One potential issue is that $R_c(\mathbf{K})$ involves the conditional expectation with respect to the past trajectory $\mathcal{H}_t$. Fortunately, it is known that the mean-variance risk metric allows for a quadratic form reformulation, as introduced in \cite{tsiamis,tsiamis2}. {\KBB Under the linearized dynamics \eqref{eq:dynamics} and a reasonable assumption that $\bm{\xi}'_t$ has a finite fourth-order moment, we have:}
\begin{align}
	R_c(\mathbf{K}) \!= \!\! \lim_{T \rightarrow \infty} \frac{1}{T} \mathbb{E} \sum_{t=0}^{T-1} \!\Big(4\mathbf{x}_t^\intercal \mathbf{Q}\mathbf{W}\mathbf{Q} \mathbf{x}_t+4\mathbf{x}_t^\intercal \mathbf{Q}\mathbf{M}_3 \Big) \!\leq\! \bar{c},\label{eq:rconst2}
\end{align}
\noindent with $\bar{c} := c-m_4 + 4 \text{tr}\{(\mathbf{W}\mathbf{Q})^2 \}$ and the (weighted) noise statistics given as
\begin{align}
	\bar{\bm{\xi}}&:=\mathbb{E}[\bm{\xi}'_t], \; \mathbf{W}:=\mathbb{E}[(\bm{\xi}'_t-\bar{\bm{\xi}})(\bm{\xi}'_t-\bar{\bm{\xi}})^\intercal], \nonumber\\ 
	\mathbf{M}_3 &:= \mathbb{E}[(\bm{\xi}'_t-\bar{\bm{\xi}})(\bm{\xi}'_t-\bar{\bm{\xi}})^\intercal \mathbf{Q} (\bm{\xi}'_t-\bar{\bm{\xi}})], \nonumber\\
	m_4&:=\mathbb{E}[(\bm{\xi}'_t-\bar{\bm{\xi}})^\intercal \mathbf{Q} (\bm{\xi}'_t-\bar{\bm{\xi}}) - \text{tr}(\mathbf{W}\mathbf{Q})]^2. \nonumber
\end{align} 
Thanks to the quadratic form in \eqref{eq:rconst2} same as \eqref{eq:obj}, the ensuing section will develop a dual approach to solve \eqref{eq:opt}. 

%

\section{Reinforcement Learning based Algorithm} \label{sec:ALG}

The risk constraint in \eqref{eq:opt} makes it challenging to find a closed-form solution. {\KBA This inspired us to develop a stochastic gradient-descent method with a max-oracle (SGDmax), which can efficiently find the stationary point (SP) of \eqref{eq:opt} while considering structured feedback.}
We first form its Lagrangian function with the multiplier $\lambda \geq 0$, as given by
\begin{align}
	&\mathcal{L}(\mathbf{K}, \lambda) 
    :=R_0(\mathbf{K})+ \lambda  [R_c(\mathbf{K})-\bar{c}]\nonumber\\
	=&\lim_{T \rightarrow \infty} \frac{1}{T}\mathbb{E}\sum_{t=0}^{T-1} \left[\mathbf{x}^\intercal_t \mathbf{Q}_{\lambda} \mathbf{x}_t\!+\!\mathbf{u}^\intercal_t \mathbf{R} \mathbf{u}_t\!+\!4 \lambda\mathbf{x}_t^\intercal \mathbf{Q} \mathbf{M}_3\right]\!-\!\lambda\bar{c}, \label{eq:lag}
\end{align}
\noindent with the matrix $\mathbf{Q}_\lambda := \mathbf{Q} + 4\lambda (\mathbf{Q} \mathbf{W} \mathbf{Q})$. The Lagrangian clearly follows the same quadratic form as the original LQR objective, thereby allowing to adopt a policy gradient approach popularly used for the unconstrained LQR problem \cite{tsiamis}. 
Based on \eqref{eq:lag}, the dual problem becomes
\begin{align} 
\max_{\lambda \in \mathcal{Y}} \mathcal{D}(\lambda) = \max_{\lambda \in \mathcal{Y}}\min_{\mathbf{K} \in \mathcal{K}} \mathcal{L}(\mathbf{K}, \lambda), \label{eq:maximin}
\end{align}

\noindent with $\mathcal{D}(\lambda)$ solved by an inner minimization problem to determine the best $\mathbf{K}$ within the structured set $\mathcal{K}$. Note that we use a bounded set $\mathcal{Y} := [0, \Lambda]$ for $\lambda$ by assuming that \eqref{eq:opt} is feasible {\KBB with reasonable bound $c$} and thus $\lambda$ is finite. 

{\KBA Solving \eqref{eq:maximin} presents a significant challenge for several reasons. First, it requires the use of a primal-dual-based gradient method, which involves both inner and outer problems. Second, we need to take into account the structured feedback when determining $\mathbf{K}\in\mathcal{K}$. Last, due to the intricate nature of \eqref{eq:lag}, the computation of first-order gradients involves a substantial computational burden.}
To resolve these issues, we consider its minimax counterpart problem, which is 
\begin{align}
    \min_{\mathbf{K} \in \mathcal{K}} \Phi(\mathbf{K}),  ~\textrm{where}~ \Phi(\mathbf{K}) := \max_{\lambda \in \mathcal{Y}} \mathcal{L}(\mathbf{K},\lambda). \label{eq:minimax}
\end{align}
Since $\mathcal{L}(\mathbf{K}, \lambda)$ is a linear function of $\lambda$ per \eqref{eq:lag}, we can directly find the optimal $\lambda$ to the inner problem in \eqref{eq:minimax}.  If the constraint is satisfied, we have $\lambda = 0$,  and otherwise, $\lambda = \Lambda$ holds. While \eqref{eq:minimax} is not the same problem as \eqref{eq:maximin}, the KKT stationary conditions for \eqref{eq:maximin} are closely related to those of \eqref{eq:minimax}. This allows us to solve \eqref{eq:minimax} instead of \eqref{eq:maximin} to find the SP of the original problem as advocated in our earlier work \cite{kwon}. 
To further minimize $\Phi(\mathbf{K})$, we can adopt the gradient descent (GD) method. However, we must take into account the structured set $\mathcal{K}$ introduced in Section~\ref{sec:RC}, which could complicate the feasible region \cite{feng}. 
Fortunately, we can compute the gradient over $\ccalK$ easily by considering only the non-zero entries, as denoted by $\nabla_{\mathcal{K}} \mathcal{L}(\mathbf{K},\lambda)$ \cite{kwon}. It is worth noting that even though $\Phi(\mathbf{K})$ is not differentiable everywhere due to the discontinuity of the optimal $\lambda$, we may still use it as the subgradient for $\Phi(\mathbf{K})$, as earlier work \cite{GDA} has established that $\nabla_{\mathcal{K}} \mathcal{L}(\mathbf{K},\lambda) \in \partial \Phi(\mathbf{K})$.

\begin{algorithm}[t]
	\SetAlgoLined
	\caption{Zero-Order Policy Gradient (ZOPG) Estimation}
	\label{alg:ZOO}
	\DontPrintSemicolon
	{\bf Inputs:} smoothing radius $r$, 
	the policy $\mathbf{K}$ and its perturbation  $\mathbf{U} \in \mathcal{S}_{\mathcal{K}}$, both have $n_\ccalK$ non-zero entries.\; 
	Obtain $\lambda' \leftarrow \argmax_{\lambda \in \ccalY} \ccalL(\mathbf{K}+r\mathbf{U}, \lambda)$;\; 
	Estimate the gradient $\hat{\nabla}_{\ccalK} \ccalL (\mathbf{K}; \mathbf{U}) = \frac{n_\ccalK}{r} \ccalL(\mathbf{K}+r\mathbf{U}, \lambda') \mathbf{U}$.\;
	{\bf Return:} $\hat{\nabla}_{\ccalK} \ccalL (\mathbf{K}; \mathbf{U})$.
\end{algorithm}
\setlength{\textfloatsep}{0.3cm}

Although it is possible to define $\nabla_{\mathcal{K}} \mathcal{L}(\mathbf{K},\lambda)$, computing this first-order term can lead to  high computational complexity. {\KBA To address this issue, we utilize the zero-order policy gradient (ZOPG) method as tabulated in Algorithm~\ref{alg:ZOO}. The ZOPG method relies on function values to provide unbiased gradient estimates without the need to consider first-order derivatives \cite{ZOO}. Since we can directly obtain function values from the trajectory, there is no requirement for additional steps to compute gradients, offering significant advantages in reducing computational workload.}
Based on the ZOPG estimate, we implement the stochastic gradient-descent with max-oracle (SGDmax) algorithm as tabulated in Algorithm~\ref{alg:SGD}. Starting from the initial policy $\mathbf{K}_0$, we perform iterative GD update on $\mathbf{K}$. To reduce estimation variance, we utilize $\hat{\mathbf{G}}(\mathbf{K})$, which represents the average of $N$ ZOPG estimates. 
Note that by determining proper smoothing radius $r$, step-size $\eta$ and the number of iterations $M$ based on the Lipschitz and smoothness constants of $\Phi(K)$, we can attain a high convergence probability to SP, approximately 90\% \cite{kwon}. This allows us to implement Algorithm~\ref{alg:SGD} conveniently in practice, as demonstrated in the ensuing section. 

\begin{algorithm}[t]
	\SetAlgoLined
	\caption{Stochastic gradient-descent with max-oracle (SGDmax)}
	\label{alg:SGD}
	\DontPrintSemicolon
	{\bf Inputs:} A feasible and stable policy $\mathbf{K}_0$, upper bound $\Lambda$ for $\lambda$, step-size $\eta$, the number of ZOPG samples $N$ and the number of iterations $M$.\; 
	\For{$m = 0, 1, \ldots, M-1$}{
		\For{$s=1,\ldots,N$}{
		Sample the random $\mathbf{U}_s \in \ccalS_\ccalK$. \;
		Use Algorithm~\ref{alg:ZOO} to return $\hat{\nabla}\ccalL_\ccalK(\mathbf{K}_m;\mathbf{U}_s)$. 
		}
		Compute the averaged stochastic gradient $\hat{\mathbf{G}}(\mathbf{K}_m) = \frac{1}{N} \sum_{s=1}^N \hat{\nabla} \ccalL(\mathbf{K}_m;\mathbf{U}_s)$. \;
		Update $\mathbf{K}_{m+1}\!\leftarrow\! \mathbf{K}_m\!-\!\eta \hat{\mathbf{G}}(\mathbf{K}_m)$.}
	{\bf Return:} the final iterate $\mathbf{K}_m$.
\end{algorithm}


%

\section{Numerical Experiments} \label{sec:NT}

\begin{figure}[t]
	\centering
	\includegraphics[width=0.9\linewidth]{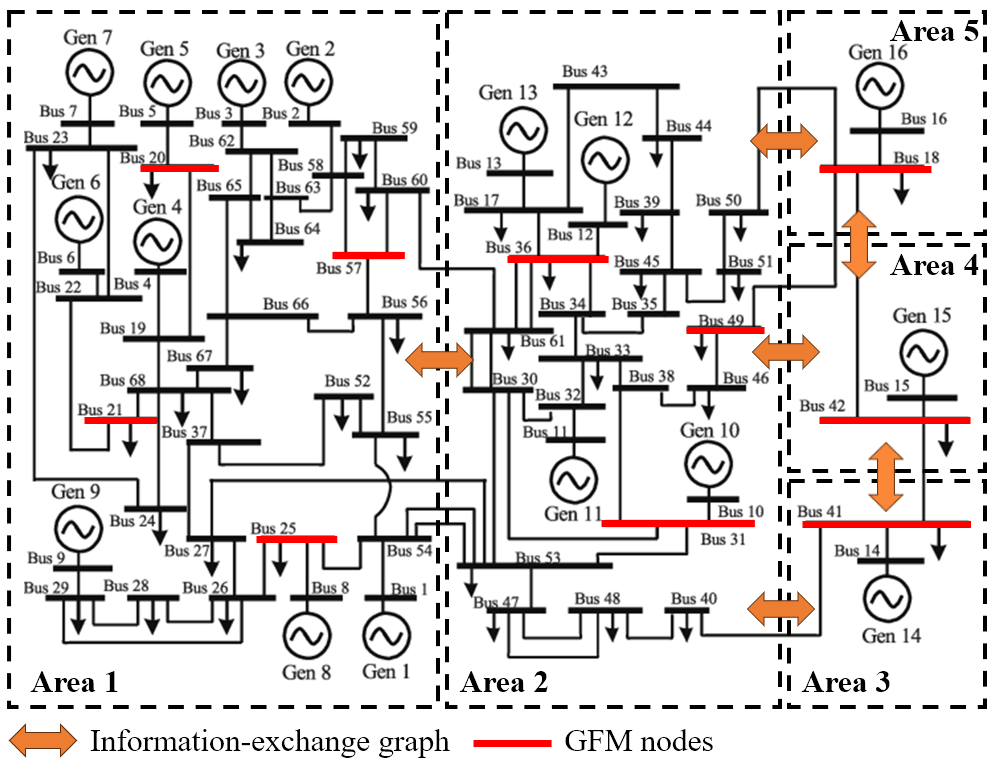}
	\caption{IEEE 68-bus system with 16 SGs and 10 GFMs.}
	\label{fig:system}
\end{figure}

To demonstrate the effectiveness of our risk-aware GFM control, we conduct numerical tests on a modified IEEE 68-bus system \cite{ieee68}. As shown in Fig.~\ref{fig:system}, this system consists of 16 SGs across five ares.
We have added a total of 10 GFMs to selected load buses. {\KB The GFM parameters follow from \cite{du} and are set to: 
$\tau = 0.01$~s, $m^p = 0.01$~pu, $m^q = 0.05$~pu, $k^{pv} = 0.01$~pu, $k^{iv} = 5.86$~pu/s.} We assume that only SGs and GFMs in neighboring areas can exchange data, forming a structured feedback following the information-exchange graph. 

We consider three cases, namely \textit{Baseline}, \textit{GFM} and \textit{GFM-Risk}. {\KB \textit{Baseline} refers to the discrete-time model-based LQR solution of \eqref{eq:obj} without taking into account {\KBA both the structured feedback and perturbations $\bm{\xi}_t$ in \eqref{eq:dynamics1}.} {\FINAL \textit{GFM} represents the typical model-based RL policy, which has been trained using Algorithm~\ref{alg:SGD} to minimize the risk-neutral \eqref{eq:obj} under \eqref{eq:dynamics1} without the risk constraint.}
\textit{GFM-Risk} is the solution of our risk-constrained problem \eqref{eq:opt} by utilizing Algorithm~\ref{alg:SGD}.} Here, \textit{GFM} and \textit{GFM-Risk} utilize Algorithm~\ref{alg:SGD} with estimated gradients obtained from Algorithm~\ref{alg:ZOO}. {\KBA The parameters used in the simulation are as follows: $r = 0.1, M = 50, \eta = 10^{-4}$ and $c = 0.2$ by considering \cite[Theorem~2]{kwon}}. The control time step is set to $\Delta t = 0.01 s$ and the observation time window is set as $0-6s$. As for perturbations, we introduce step load changes in the GFM buses at $t=0$, ranging from $\pm 0.5$~pu and $\pm 1$~pu, which we refer to as the low and high perturbation settings.

\begin{figure}[t]
	\centering
	\includegraphics[width=0.9\linewidth]{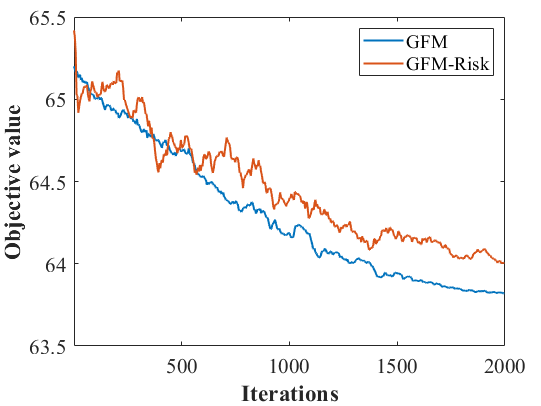}
	\caption{Trajectory of the training objective values for GFM and GFM-Risk.}
	\label{fig:train}
\end{figure}

\textbf{Training comparisons:} Fig.~\ref{fig:train} illustrates the objective trajectories for \textit{GFM} and \textit{GFM-Risk}. We can observe that both cases converge to SP using the proposed algorithm, but \textit{GFM-Risk} exhibits more fluctuations in its trajectory compared to \textit{GFM}. {\KBA This is because considering risk constraint complicates the feasible region, which leads to oscillations in $\lambda$. As $\lambda$ fluctuates, the function value $\mathcal{L}(\mathbf{K}, \lambda)$ also oscillates, consequently affecting the ZOPG that is proportional to this value.} {\KBA \textit{Baseline} is not shown in the figure, as it does not require a training.}

\begin{figure}[t]
    \centering
    \subfigure[]{\label{fig:freq_low}\includegraphics[width=0.9\linewidth]{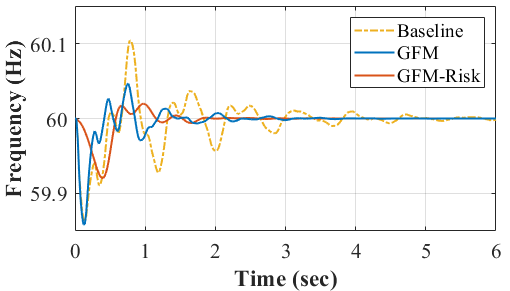}}
    \subfigure[]{\label{fig:freq_high}\includegraphics[width=0.9\linewidth]{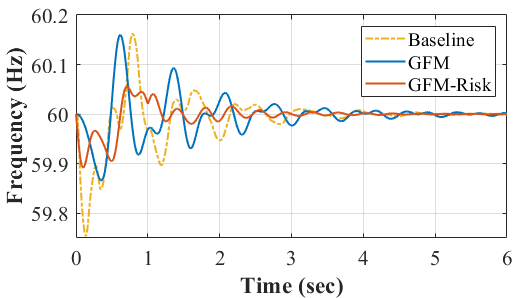}}
    \caption{Comparison on the frequency deviation of bus 67 in MG1 at (a) low and (b) high perturbation settings.}
    \captionsetup{justification=centering}
    \label{fig:freq}
\end{figure}

\textbf{Testing comparisons:} Using the three different policies from \textit{Baseline}, \textit{GFM} and \textit{GFM-Risk}, we conduct tests with 100 new scenarios, each involving random load changes of either low or high perturbation settings. To show the effectiveness of the risk constraint under extreme load changes, we select one specific scenario with the most significant changes among 100 scenarios for each setting. Fig.~\ref{fig:freq} compares the frequency deviations of bus 57 in Area 1 with three different polices in the both perturbation settings.

\begin{figure}[t]
    \centering
    \includegraphics[width=\linewidth]{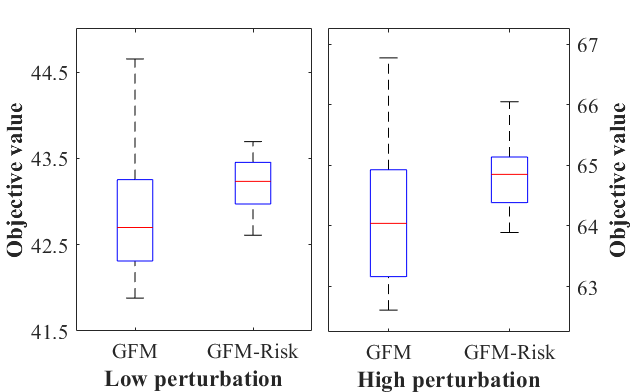}
    \caption{Comparison on the LQR objective values in low and high load perturbations.}
    \captionsetup{justification=centering}
    \label{fig:boxplot}
\end{figure}

Clearly, both \textit{GFM} and \textit{GFM-Risk} effectively mitigate frequency deviations and reach a steady-state faster in all settings compared to \textit{Baseline}. Here, Fig~\ref{fig:freq_high} experiences more fluctuations than Fig~\ref{fig:freq_low} in both cases, thus taking more time to reach the steady-state. This is reasonable since the maximum perturbation level increases from $\pm 0.5$~pu to $\pm 1$~pu. In addition, \textit{GFM-Risk} demonstrates smaller fluctuations and leads to the steady-state faster than \textit{GFM} in both settings. This implies that \textit{GFM-Risk} provides more damping in extreme cases compared to \textit{GFM} due to the risk constraint. Notably, we observe an improvement in the damping performance of \textit{GFM-Risk} in Fig.~\ref{fig:freq_high}, indicating that the constraint becomes more effective as perturbation level increases.


To further analyze the quantitative results of implementing the risk constraint, we present Fig.~\ref{fig:boxplot}, which provides statistical information regarding the LQR objective values obtained from 100 testing scenarios. The red lines indicate the median values while the lower and upper quartiles are represented by the blue boxes. The maximum and minimum values are depicted with black lines. 
First, it is evident that both \textit{GFM} and \textit{GFM-Risk} in low load perturbation exhibit larger median values and variances compared to high perturbation. This is expected since an increase in the maximum perturbation levels introduces higher variability into the system, resulting in more frequency oscillations. Second, \textit{GFM-Risk} exhibits a slightly higher median value compared to \textit{GFM} in both settings, as expected due to the constraint. However, \textit{GFM-Risk} demonstrates a significantly smaller variance, notably reducing the gap between the upper quartile and the maximum value. As a higher objective value implies more frequency deviation, this result showcases the effectiveness of incorporating the risk constraint in mitigating worst-case performance and enhancing system stability. Last, the variance of \textit{GFM-Risk} decreases dramatically in a high perturbation setting, validating that the constraint becomes more effective in environments with significant load changes. 
%
%
%
%
\section{Conclusions} \label{sec:CON}
This paper designed a risk-aware controller for GFMs that aims to address frequency oscillations resulting from high load perturbations in GFM inverter-dominated power systems. Based on the linearized system model that incorporates both SGs and GFMs, we formulated the problem to minimize the LQR objective over the structured feedback gain according to the connectivity of communication network. Since increased load perturbations lead to the higher state variability, we introduced the mean-variance risk constraint to limit the state cost variations, thereby reducing the system variability and enhancing the worst-case performance. By reformulating this constraint into a tractable quadratic form, we solved the dual problem using an RL-based SGDmax algorithm. This method searched for a policy through GD iterations by leveraging efficient ZOPG for gradient estimation. Numerical tests on the modified IEEE 68-bus system highlighted the effectiveness of the proposed risk-aware GFM controller in reducing the variability of total LQR cost, thus improving the performance in worst-case scenarios with significant load perturbations.


Future directions include exploring additional risk measures such as conditional value at risk (CVaR), as well as examining
the impact of other types of perturbations. 
{\KBA In addition, it is worth highlighting that our risk-aware control framework has a variety of applications in other dynamic systems in mitigating the system variability in practice.} 
 
\bibliographystyle{IEEEtran}
\bibliography{bibliography}

\begin{thebibliography}{10}
\providecommand{\url}[1]{#1}
\csname url@samestyle\endcsname
\providecommand{\newblock}{\relax}
\providecommand{\bibinfo}[2]{#2}
\providecommand{\BIBentrySTDinterwordspacing}{\spaceskip=0pt\relax}
\providecommand{\BIBentryALTinterwordstretchfactor}{4}
\providecommand{\BIBentryALTinterwordspacing}{\spaceskip=\fontdimen2\font plus
\BIBentryALTinterwordstretchfactor\fontdimen3\font minus
  \fontdimen4\font\relax}
\providecommand{\BIBforeignlanguage}[2]{{%
\expandafter\ifx\csname l@#1\endcsname\relax
\typeout{** WARNING: IEEEtran.bst: No hyphenation pattern has been}%
\typeout{** loaded for the language `#1'. Using the pattern for}%
\typeout{** the default language instead.}%
\else
\language=\csname l@#1\endcsname
\fi
#2}}
\providecommand{\BIBdecl}{\relax}
\BIBdecl

\bibitem{anttila}
S.~Anttila, J.~S. Döhler, J.~G. Oliveira, and C.~Boström, ``{Grid Forming
  Inverters: A Review of the State of the Art of Key Elements for Microgrid
  Operation},'' \emph{Energies}, vol.~15, no.~15, 2022.

\bibitem{long}
A.~Singhal, T.~L. Vu, and W.~Du, ``{Consensus Control for Coordinating
  Grid-Forming and Grid-Following Inverters in Microgrids},'' \emph{IEEE
  Transactions on Smart Grid}, vol.~13, no.~5, pp. 4123--4133, 2022.

\bibitem{rath}
D.~B. Rathnayake, M.~Akrami, C.~Phurailatpam, S.~P. Me, S.~Hadavi,
  G.~Jayasinghe, S.~Zabihi, and B.~Bahrani, ``{Grid Forming Inverter Modeling,
  Control, and Applications},'' \emph{IEEE Access}, vol.~9, pp.
  114\,781--114\,807, 2021.

\bibitem{wei}
W.~Du, Z.~Chen, K.~P. Schneider, R.~H. Lasseter, S.~P. Nandanoori, F.~K.
  Tuffner, and S.~Kundu, ``{A Comparative Study of Two Widely Used Grid-forming
  Droop Controls on Microgrid Small-Signal Stability},'' \emph{IEEE Journal of
  Emerging and Selected Topics in Power Electronics}, vol.~8, no.~2, pp.
  963--975, 2019.

\bibitem{peyghami}
S.~Peyghami, P.~Davari, H.~Mokhtari, and F.~Blaabjerg, ``{Decentralized Droop
  Control in DC Microgrids based on a Frequency Injection Approach},''
  \emph{IEEE Transactions on Smart Grid}, vol.~10, no.~6, pp. 6782--6791, 2019.

\bibitem{ebrahimi}
M.~Ebrahimi, S.~A. Khajehoddin, and M.~Karimi-Ghartemani, ``{An Improved
  Damping Method for Virtual Synchronous Machines},'' \emph{IEEE Transactions
  on Sustainable Energy}, vol.~10, no.~3, pp. 1491--1500, 2019.

\bibitem{serban}
I.~Serban and C.~P. Ion, ``{Microgrid control based on a grid-forming inverter
  operating as virtual synchronous generator with enhanced dynamic response
  capability},'' \emph{International Journal of Electrical Power \& Energy
  Systems}, vol.~89, pp. 94--105, 2017.

\bibitem{grob}
D.~Groß, M.~Colombino, J.-S. Brouillon, and F.~Dörfler, ``{The Effect of
  Transmission-Line Dynamics on Grid-Forming Dispatchable Virtual Oscillator
  Control},'' \emph{IEEE Transactions on Control of Network Systems}, vol.~6,
  no.~3, pp. 1148--1160, 2019.

\bibitem{du}
W.~Du, Y.~Liu, R.~Huang, F.~K. Tuffner, J.~Xie, and Z.~Huang,
  ``{Positive-Sequence Phasor Modeling of Droop-Controlled, Grid-Forming
  Inverters with Fault Current Limiting Function},'' in \emph{{IEEE Power \&
  Energy Society Innovative Smart Grid Technologies Conference}}, 2022, pp.
  1--5.

\bibitem{behrooz}
B.~Mirafzal and A.~Adib, ``{On Grid-Interactive Smart Inverters: Features and
  Advancements},'' \emph{IEEE Access}, vol.~8, pp. 160\,526--160\,536, 2020.

\bibitem{li}
Y.~Li, W.~Gao, S.~Huang, R.~Wang, W.~Yan, V.~Gevorgian, and D.~W. Gao,
  ``{Data-Driven Optimal Control Strategy for Virtual Synchronous Generator via
  Deep Reinforcement Learning Approach},'' \emph{Journal of Modern Power
  Systems and Clean Energy}, vol.~9, no.~4, pp. 919--929, 2021.

\bibitem{mukherjee2021scalable}
S.~Mukherjee, A.~Chakrabortty, H.~Bai, A.~Darvishi, and B.~Fardanesh,
  ``{Scalable designs for reinforcement learning-based wide-area damping
  control},'' \emph{IEEE Transactions on Smart Grid}, vol.~12, no.~3, pp.
  2389--2401, 2021.

\bibitem{lintao}
L.~Ye, H.~Zhu, and V.~Gupta, ``{On the Sample Complexity of Decentralized
  Linear Quadratic Regulator With Partially Nested Information Structure},''
  \emph{IEEE Transactions on Automatic Control}, vol.~68, no.~8, pp.
  4841--4856, 2023.

\bibitem{mbrl1}
H.~Issa, V.~Debusschere, L.~Garbuio, P.~Lalanda, and N.~Hadjsaid, ``Artificial
  intelligence-based controller for grid-forming inverter-based generators,''
  in \emph{2022 IEEE PES Innovative Smart Grid Technologies Conference Europe
  (ISGT-Europe)}, 2022, pp. 1--6.

\bibitem{mbrl2}
C.~Zhao, C.~Wang, Z.~Cheng, W.~Sun, and B.~Zhao, ``Optimal distributed
  coordinated reinforcement learning for secondary voltage control in
  time-delayed microgrid,'' \emph{IEEE Systems Journal}, vol.~17, no.~3, pp.
  3480--3491, 2023.

\bibitem{eskandari}
M.~Eskandari, A.~V. Savkin, and J.~Fletcher, ``{A Deep Reinforcement
  Learning-based Intelligent Grid-Forming Inverter for Inertia Synthesis by
  Impedance Emulation},'' \emph{IEEE Transactions on Power Systems}, 2023.

\bibitem{kwon}
K.-b. Kwon, L.~Ye, V.~Gupta, and H.~Zhu, ``Model-free learning for
  risk-constrained linear quadratic regulator with structured feedback in
  networked systems,'' in \emph{IEEE 61st Conference on Decision and Control
  (CDC)}, 2022, pp. 7260--7265.

\bibitem{kundur}
P.~Kundur, \emph{Power System Stability and Control}.\hskip 1em plus 0.5em
  minus 0.4em\relax McGraw-Hill, 1994.

\bibitem{kron}
V.~Vittal, J.~D. McCalley, P.~M. Anderson, and A.~Fouad, \emph{{Power System
  Control and Stability}}.\hskip 1em plus 0.5em minus 0.4em\relax John Wiley \&
  Sons, 2019.

\bibitem{feng}
H.~Feng and J.~Lavaei, ``{On the Exponential Number of Connected Components for
  the Feasible Set of Optimal Decentralized Control Problems},'' in \emph{2019
  American Control Conference (ACC)}, 2019, pp. 1430--1437.

\bibitem{tsiamis}
A.~Tsiamis, D.~S. Kalogerias, L.~F.~O. Chamon, A.~Ribeiro, and G.~J. Pappas,
  ``{Risk-Constrained Linear-Quadratic Regulators},'' in \emph{{59th IEEE
  Conference on Decision and Control (CDC)}}, 2020, pp. 3040--3047.

\bibitem{tsiamis2}
A.~Tsiamis, D.~S. Kalogerias, A.~Ribeiro, and G.~J. Pappas, ``{Linear Quadratic
  Control with Risk Constraints},'' \emph{arXiv preprint arXiv:2112.07564},
  2021.

\bibitem{GDA}
T.~Lin, C.~Jin, and M.~Jordan, ``{On Gradient Descent Ascent for
  Nonconvex-Concave Minimax Problems},'' in \emph{International Conference on
  Machine Learning}.\hskip 1em plus 0.5em minus 0.4em\relax PMLR, 2020, pp.
  6083--6093.

\bibitem{ZOO}
J.~C. Spall, ``{A One-Measurement Form of Simultaneous Perturbation Stochastic
  Approximation},'' \emph{Automatica}, vol.~33, no.~1, pp. 109--112, 1997.

\bibitem{ieee68}
B.~Pal and B.~Chaudhuri, \emph{{Robust Control in Power Systems}}.\hskip 1em
  plus 0.5em minus 0.4em\relax Springer Science \& Business Media, 2006.

\end{thebibliography}

\end{document}